  \documentclass[reqno]{amsart}   
\usepackage{pb-diagram,enumerate,amsmath}   
\usepackage{pstricks,pst-node,pst-coil}          
           
\def\al{{\alpha}}           
\def\be{{\beta}}           
\def\de{{\delta}}           
           
\def\om{{\omega}}

\def\Si{{\Sigma}}           
\def\ga{{\gamma}}           
\def\ep{{\varepsilon}}           
\def\Ga{{\Gamma}}

\def\Th{{\Theta}}           
           
\def\phi{{\varphi}}           
\def\Ph{{\Phi}}

\def\na{{\nabla}}           
\def\pa{{\partial}}

\def\proof#1{{\par\medbreak\noindent {\bf Proof\setbox0\hbox{#1}%
\ifdim\wd0=0pt .\else\ \ignorespaces #1.\fi}\enspace}}           
\def\iop#1{{\par\medbreak\noindent {\bf Idea of proof\setbox0\hbox{#1}%
\ifdim\wd0=0pt .\else\ \ignorespaces #1.\fi}\enspace}}           
\def\examples{{\noindent {\bf Examples. }\par\kern-\baselineskip}}           

\DeclareMathAlphabet{\doba}{U}{msb}{m}{n}

\gdef\mR{\doba{R}}

\def\qed{{\leavevmode\unskip\nobreak\hfil\penalty 50\hskip 1em%
  \hbox{}\nobreak\hfil\lower 1pt\hbox{$\Box$\kern-.5pt}\parfillskip 0pt           
  \finalhyphendemerits 0\par\bigbreak}}           
\def\qedmath#1{\setbox0\hbox{$\displaystyle #1$}\templaenge=\textwidth\advance\templaenge by -\wd0%
\setbox1\hbox{$\Box$}\advance\templaenge by -2\wd1%
$$#1\hbox to0pt{\kern.5\templaenge$\Box$\kern-.5pt\hss}$$\par\bigbreak}

\long\def\komment#1{}           
\newtheoremstyle{remarks}{3pt}{3pt}{}{}{\bfseries}{}{ }{}           
           
\makeatletter           
\@addtoreset{equation}{section}           
\renewcommand{\theequation}{\the\c@section.\the\c@equation}           
\makeatother           
\newtheorem{theorem}[equation]{\bf T{\footnotesize HEOREM}}           
\newtheorem{proposition}[equation]{\bf P{\footnotesize ROPOSITION}}           
\newtheorem{lemma}[equation]{\bf L{\footnotesize EMMA}}

\theoremstyle{definition}           
\newtheorem*{remark}{Remark}           
           
\newtheorem*{definition}{Definition}

\theoremstyle{remarks}           

\let\<\langle           
\let\>\rangle

\def\vol{{\mathop{{\rm vol}}}}           
           
\def\scal{{\mathop{{\rm Scal}}}}           
           
\def\Ric{{\mathop{{\rm Ric}}}}           
\let\witi\widetilde           
\def\can{{\mathop{\rm can}}}           
           
\begin{document}           
\title{Positive mass theorem for the Yamabe problem on spin manifolds}           
\author{Bernd Ammann, Emmanuel Humbert}           
\date{Dec 15, 2005}           
\maketitle           
              
\begin{abstract}         
Let $(M,g)$ be a compact connected spin manifold of dimension          
$n\geq 3$ whose Yamabe invariant is positive. We assume          
that  $(M,g)$ is locally conformally flat or that $n \in         
\{3,4,5\}$. According to a positive mass theorem  
by Schoen and Yau  
the constant term in the asymptotic          
development of the Green's function of         
the conformal Laplacian is positive if $(M,g)$ is not conformally         
equivalent to the sphere.  The proof was simplified by Witten with  
the help of spinors. 
In our article we will give a proof which is even   
considerably shorter.   
Our proof is a modification of Witten's argument, but no analysis    
on asymptotically flat spaces is needed.   
\end{abstract}           
           
           
{\bf Mathematics Classification:}           
53C21 (Primary), 58E11, 53C27 (Secondary)           
           
\section{Introduction}           
           
The positive mass conjecture is a famous and difficult problem which  
originated in physics.           
The mass is a Riemannian invariant of an asymptotically           
flat manifold of dimension $n \geq 3$ and of order            
$\tau > {n-2 \over2}$.           
The problem consists in proving that the mass is             
positive if the manifold is not conformally diffeomorphic to           
$(\mR^n,\hbox{can})$.            
Two good references on this subject are \cite{lee.parker:87,herzlich:98}.       
           
Schoen and Yau \cite{schoen:89,schoen.yau:79a}            
gave a proof if the dimension is at most $7$           
and Witten \cite{witten:81,bartnik:86}            
proved the result if the manifold is spin.            
The positivity of the mass has been proved in several other            
particular cases (see e.g. \cite{schoen:84}),           
but the conjecture in its full generality still remains open.                 
           
This problem played an important role in geometry because its solution  
led to the solution of the Yamabe problem.            
Namely, let $(M,g)$ be a compact connected Riemannian manifold of dimension            
$n\geq 3$.            
In \cite{yamabe:60} Yamabe attempted to show that there is a metric           
$\tilde{g}$           
 conformal to $g$ such that the scalar curvature            
$Scal_{\tilde{g}}$ of ${\tilde{g}}$ is            
constant. However, Trudinger realized that Yamabe's            
proof contained a serious gap.           
It was the achievement of many mathematicians to finally solve the problem            
of finding a conformal metric ${\tilde{g}}$ with constant scalar curvature.           
The problem of finding a conformal ${\tilde{g}}$ with constant scalar curvature           
is called the Yamabe problem.              
\noindent As a first step, Trudinger \cite{trudinger:68}           
was able to repair the gap if a conformal invariant named the Yamabe           
invariant is non-positive.  
The problem is much more difficult if the Yamabe invariant is positive,            
which is equivalent to the existence of a metric of positive            
scalar curvature in the conformal class of $g$.            
Aubin \cite{aubin:76} solved the problem  when  $n \geq6$ and  $M$ is not locally conformally flat.            
Then, in \cite{schoen:84}, Schoen completed the proof that a solution to the  
Yamabe problem exists by            
using the positive mass theorem in the remaining           
cases.          
Namely, assume           
that $(M,g)$ is locally conformally flat or $n \in \{ 3,4,5\}$.           
Let            
$$L_g= {4(n-1) \over n-2} \Delta_g + Scal_g$$           
be the conformal Laplacian of the metric $g$ and $P \in M$.             
There exists a smooth function $\Gamma$, the so-called Green's           
function of $L_g$, which is defined on $M -\{P \}$            
such that           
$L_g\Gamma = \delta_P$ in the sense of distributions (see for example           
\cite{lee.parker:87} ).            
Moreover, if we let $r=d_g(.,P)$, then in conformal normal coordinates    
$\Gamma$ has the following expansion           
at $P$:           
  $$ \Gamma(x)=        
    {1 \over 4(n-1)\om_{n-1}\,r^{n-2} } +A + \alpha(x)\qquad \om_{n-1}=\vol(S^{n-1})$$           
where $A \in \mR$.            
In addition,  $\alpha$ is a function defined on a           
neighborhood of $P$ and $\alpha(0)=0$. On this neighborhood of $P$, the function  
$\alpha$ is smooth         
if $(M,g)$ is locally conformally           
flat, and it is a Lipschitz function for $n=3,4,5$.           
Hence, in both cases $\alpha=O(r)$.            
Schoen has shown in \cite{schoen:84} that the positivity of $A$           
would imply the solution of the Yamabe problem. He also proved that            
$A$  is a positive multiple of the mass of the asymptotically flat manifold            
$(M,\Gamma^{{4 \over n-2}}g)$.            
Hence, in these special cases the solution of the Yamabe problem follows from the positive mass       
theorem, which was proven by Schoen and Yau            
in \cite{schoen.yau:79a,schoen.yau:88}.           
\\           
           
In our article, we will give a short proof for the positivity            
of the constant term $A$ in the development of the Green's function            
in case that $M$ is spin and locally conformally flat.            
The statement of this paper is weaker than the results of Witten           
\cite{witten:81} and Schoen and Yau \cite{schoen:89,schoen.yau:79a}.           
The proof in our paper is inspired by Witten's reasoning,           
but we have considerably simplified many of the analytic arguments.           
Witten's argument is based on the construction of a test spinor on the            
stereographic blowup which is both harmonic and            
asymptotically constant. We show that the Green's function           
for the Dirac operator on $M$ can be used to construct such a test spinor.            
In this way, we obtain a very short            
solution of the Yamabe problem using only elementary and well known            
facts from analysis on compact manifolds.          
         
The last section shows how to adapt our proof to arbitrary spin manifolds         
of dimensions $3$, $4$ and $5$. In dimension $3$ the          
proof is completely analogous. However, in dimensions $4$ and $5$, additional          
estimates have to be derived in order to get sufficient          
control on the Green's function of the Dirac operator.         

{\it Remark about this version:} In printed version that appeared in 
Geom. Funct. Anal., {\bf 15}, {\it 567--576} (2005) a term in the local formula
for the Dirac operator was missing. This implies that some additional 
terms have to be added in the last proof. This gap is repaired in the 
present version.

\section{The locally conformally flat case}           
           
In this section, we will assume that $(M,g)$ is a compact, connected,           
locally conformally flat spin manifold of dimension $n\geq 3$.            
The Dirac operator is denoted by~$D$. A spinor $\psi$ is called            
$D$-\emph{harmonic} if $D\psi\equiv 0$.          
As the solution of the Yamabe problem          
in the case of non-negative Yamabe invariant         
follows from \cite{trudinger:68} we will assume that the Yamabe invariant is positive.          
Hence, the conformal class contains a positive scalar curvature metric. As $\dim\ker D$ is conformally invariant,         
we see that $\dim \ker D=0$.           
We fix a point $P\in M$. We can assume that $g$ is flat in a small ball           
$B_P(\delta)$ of radius $\de$ about $P$, and that $\de$ is smaller   
than the injectivity radius. Let $(x^1,\ldots,x^n)$            
denote local coordinates on $B_P(\delta)$.           
On $B_P(\delta)$ we trivialize the spinor bundle via parallel transport.           
           
\begin{lemma}  \label{lem}         
Let $\psi_0\in \Si_P M$. Then there is a $D$-harmonic spinor $\psi$ on $M\setminus \{P\}$ satisfying           
  $$\psi|_{B_P(\de)}= \frac{x}{r^{n}}           
 \cdot \psi_0 + \theta(x) $$           
where $\theta(x)$ is a smooth spinor on $B_P(\delta)$.           
\end{lemma}           
           
\noindent It is not hard to see, that    
in the sense of distributions    
  $$D\psi= -\om_{n-1}\de_P \psi_0,$$   
where  $\de_P$ is the $\delta$-function centered at $P$.   
Hence, by definition, $-\om_{n-1}^{-1}\psi$ is           
the Green's function of the Dirac operator.           
           
\proof{}           
Our construction of $\psi$ follows the construction of the Green's function            
$G$ of the Laplacian in \cite[Lemma~6.4]{lee.parker:87}.            
Namely, we take a cut-off function $\eta$ with support in $B_P(\de)$           
which is equal to $1$ on $B_P(\de/2)$.    
We set $\Ph= \eta \frac{1}{r^{n-1}}           
\frac{ x} {r} \cdot \psi_0$ where $\psi_0$ is constant. The spinor           
$\Ph$ is $\bar{D}$-harmonic on $B_P(\de/2)\setminus \{P\}$.         
Outside $B_P(\delta)$ we extend $\Ph$ by zero, and we obtain a 
smooth spinor on $M\setminus \{P\}$.         
As $\bar{D}\Phi|_{B_P(\de/2)}\equiv 0$, we see that $\bar{D}\Phi$ extends to        
a smooth spinor on $M$.         
Using the selfadjointness of $\bar{D}$ together with $\ker \bar{D}=\{0\}$
we know that           
there is a smooth spinor $\theta_1$ such that    
$\bar{D} \theta_1= -D\Ph$.    
Obviously, $\psi= \Ph +\theta_1$ is a spinor as claimed.           
\qed           
       
       
We now show that the existence of $\psi$ implies            
the positivity of $A$.

\begin{theorem} \label{mass}           
Let $(M,g)$ be a compact connected locally conformally flat manifold of 
dimension $n \geq 3$. Then, the mass $A$ of $(M,g)$ satisfies $A \geq 0$.  
Furthermore, equality holds            
if and only if $(M,g)$           
is conformally equivalent to the standard sphere $(S^n, \can)$.           
\end{theorem}             
           
\proof{}             
Let $\psi$ be given by lemma \ref{lem}. Without loss of generality, we   
may assume that $|\psi_0|=1$. Let $\Gamma$ be the Green's function for   
$L_g$,   
and $G = 4(n-1)\om_{n-1}\Gamma$.       
Using the maximum principle, it is easy to see that        
$G$ is positive          
\cite[Lemma~6.1]{lee.parker:87}.        
We set            
$$\witi{g} = G^{\frac{4}{n-2}} g.$$            
Using the transformation formula for        
$\scal$ under conformal changes, we obtain $\scal_{\witi{g}}=0$.           
We can identify spinors on ($M\setminus \{P\}$,$\witi g$)            
with spinors ($M\setminus \{P\}$,$g$) such that         
the fiber wise scalar product on spinors is preserved \cite{hitchin:74,hijazi:86}.         
Because of the formula for the conformal change of Dirac operators,           
the spinor            
$$\witi\psi:= G^{-\frac{n-1}{n-2}}\psi$$           
is a $D$-harmonic spinor on ($M\setminus \{P\}$,$\witi g$), i.e.\ if we write           
$\witi{D}$ for the Dirac operator in the metric           
$\witi{g}$, we have $\witi{D} \witi\psi=0$.           
By the Schr\"odinger-Lichnerowicz formula we have           
 $$0= {\witi D}^2\witi\psi= \witi\na^* \witi\na \witi\psi +            
   {\scal_{\witi g}\over 4} \witi\psi = \witi\na^* \witi\na \witi\psi.$$            
Integration over $M\setminus B_P(\epsilon)$, $\epsilon >0$ and            
integration by parts yields           
$$ 0 = \int_{M\setminus B_P(\epsilon)} (\witi\na^* \witi\na \witi\psi, \witi\psi)  dv_{\witi{g}}           
     = \int_{M\setminus B_P(\epsilon)} | \witi\na \witi\psi |^2 dv_{\witi{g}}           
-\int_{S_P(\epsilon)} (\witi\na_{\witi\nu} \witi\psi,\witi\psi)           
       ds_{\witi{g}}.$$           
Here $S_P(\epsilon)$ denotes the boundary $\pa B_P(\epsilon)$,           
$\witi\nu $ is the unit normal vector on $S_P(\epsilon)$            
with respect to $\witi g$ pointing into the ball,           
and $ds_{\witi{g}}$ denotes the Riemannian volume            
element of $S_P(\epsilon)$.            
Hence, we have proved that            
\begin{eqnarray} \label{ine}           
  \int_{M\setminus B_P(\epsilon)} | \witi\na \witi\psi |^2 dv_{\witi{g}}=  \frac12  \int_{S_P(\epsilon)} \pa_{\witi\nu}  | \witi\psi           
 |^2              
ds_{\witi{g}}           
\end{eqnarray}           
If $\ep$ is sufficiently small, we have           
\begin{eqnarray} \label{normvec}           
\witi\nu=-G^{-\frac{2}{n-2}} \frac{\pa}{\pa r}=           
- \big( \ep^{2} + o(\ep^{2}) \big)           
\frac{\pa}{\pa r}            
\end{eqnarray}           
\begin{eqnarray} \label{vol}           
ds_{\witi{g}}= G^{\frac{2(n-1)}{n-2}}            
ds_g=  G^{\frac{2(n-1)}{n-2}} \ep^{n-1} ds =           
\big(\ep^{-(n-1)}+o(\ep^{-(n-1)} ) \big) ds           
\end{eqnarray}           
where $ds$ stands for the volume element of $(S^{n-1},\can)$, and            
$$| \witi\psi|^2= G^{-2\frac{n-1}{n-2}} |\psi|^2            
= \left( {1 \over r^{n-2}} +4 (n-1)\om_{n-1}\,A + r \al_1(r) \right)^{-2\frac{n-1}{n-2}}           
\left| \frac{1}{r^{n-1}} \frac{ x}{r} \cdot \psi_0+ \theta(x)           
\right|^2$$           
where $\al_1$ is a smooth function.           
This gives           
$$ | \witi\psi|^2=             
(1+4 (n-1)\om_{n-1}\,Ar^{n-2} + r^{n-1} \al_1(r) )^{-2\frac{n-1}{n-2}}           
\times$$           
$$\left( 1 + 2r^{n-1} Re(<\frac{x}{r}\cdot \psi_0, \theta(x) >) + r^{2(n-1)}            
| \theta(x)  |^2 \right)$$           
           
\noindent Noting that          
$\na_r (\frac{x}{r}\psi_0)=0$,         
we get that on               
$S_P(\epsilon)$ and  for $\ep$ small,           
\begin{eqnarray} \label{deriv}           
\pa_r  | \witi\psi|^2=  - 8 (n-1)^2\om_{n-1}\, A \ep^{n-3} + o(\ep^{n-3})            
\end{eqnarray}           
Plugging (\ref{normvec}),              
(\ref{vol}) and (\ref{deriv})           
into  (\ref{ine}),            
we get that for $\ep$ small            
\begin{eqnarray} \label{final}            
0 \leq \int_{M\setminus B_P(\epsilon)} | \witi\na \witi\psi |^2           
dv_{\witi{g}}           
 & = & 4(n-1)^2\om_{n-1}\,A  \int_{S^{n-1}} ds+ o(1)\\& =& 4 (n-1)^2            
\om_{n-1}^2 A            
+ o(1)           
\end{eqnarray}           
This implies that $A \geq 0$.            
           
Now, we assume that $A = 0$. Then it           
follows from (\ref{final})  that            
$ \witi\na \witi\psi=0$ on $M\setminus \{P\}$ and hence, $\witi\psi$ is parallel. Since the           
choice of $\psi_0$ is arbitrary, we obtain in this way a basis of           
parallel spinors on $(M\setminus \{P\}, \witi{g})$. This implies that            
$( M\setminus \{P\}, \witi{g})$ is flat and hence isometric to euclidean space.           
Let $I:(M\setminus \{P\},\witi{g})\to (\mR^n,\can)$ be an isometry. We define            
$f(x)= {1 + \|I(x)\|^2/4}, \quad x\in M$.            
Then $M\setminus \{P\}, f^{-2} \witi{g}= f^{-2}G^{4\over n-2}g$ is isometric            
to $(S^n\setminus\{N\},\can)$. The function $f^{-2}G^{4\over n-2}$ is smooth           
on $M\setminus \{P\}$ and can be extended continuously to a            
positive function on $M$. Hence, $M$ is conformal to $(S^n,\can)$.           
         
\section{The case of dimensions 3, 4 and 5} \label{345}        
Now under the assumption that the dimension of M is $3$, $4$ or $5$ we show how to adapt the proof from the last section to the case in which       
$M$ is not conformally flat.              
Let us assume that $(M,g)$ is an arbitrary connected       
spin manifold of dimension $n \in \{3,4,5\}$.         
We choose any $P \in M$. After possibly replacing $g$ by a metric    
conformal to $g$, we may assume that $\Ric_g(P)=0$.         
We trivialize the spinor bundle near $P$ with the Bourguignon-Gauduchon
trivialization \cite{bourguignon.gauduchon:92}:
let $(x_1,\cdots,x_n)$ be a system of normal coordinates at $P$ defined on
a neighborhood $V$ of $P$. Let also    
\begin{eqnarray*}   
 G: V&\longrightarrow&S^2_+(n,\mR)\\   
m&\longmapsto& G_m:=(g_{ij}(m))_{ij}    
\end{eqnarray*}   
denote the smooth map which associates to any point $m\in V$, the matrix of   
the coefficients of the metric~$g$ at this point, expressed in the basis   
$(\partial_i:=\frac{\partial}{\partial x^i})_{1\leq i\leq n}\;$.  
The vector fields $\pa_i$ are defined on a neighborhood $U$ of $0$ in
$\mR^n$. 
Since $G_m$ is    
symmetric and positive definite, there is a unique symmetric matrix  
$B_m$ such that $$B_m^2=G_m^{-1}\;.$$  

\noindent We now define  
  $$e_i:= b_i^j\partial_j\;,$$  
so that $(e_1,\ldots,e_n)$ is an   
orthonormal frame of $(TV,g)$.  
Standard constructions then allow  to  identify the spinor bundles $\Sigma
U$ and $\Sigma V$.
Denote by $\nabla$ (resp. $\bar\nabla$) the Levi-Civita connection  
on $(TU,\xi)$ (resp. $(TM,g)$) as well as its lifting    
to the spinor bundle $\Sigma U$ (resp. $\Sigma V$).  
We denote Clifford multiplication on $\Sigma V$ by ``$\;\cdot\;$''. For all    
spinor field $\psi\in\Gamma(\Sigma U)$, since   
$\bar\psi\in\Gamma(\Sigma V)$ and by definition of $\bar\nabla$, we have   
\begin{equation} 
\label{relatnabla} 
\bar\nabla_{e_i}\bar\psi=\overline{\nabla_{e_i}(\psi)}+\frac{1}{4}\sum_{j,k}   
\widetilde\Gamma_{ij}^k\,e_j\cdot   
e_k\cdot\bar\psi\;. 
\end{equation}   
where the  Christoffel symbols of the second kind  $\widetilde\Gamma^k_{ij}$ 
are defined by    
  $$\widetilde{\Gamma}^k_{ij}:=-\langle\bar\nabla_{e_i}e_j,e_k\rangle\:,$$ 

\noindent Taking Clifford multiplication by $e_i$ on each member of \eqref{relatnabla} and   
summing over $i$ yields   
\begin{equation*}   
\bar D\bar\psi=\sum_i e_i \cdot \overline{\nabla_{e_i} \psi}+
\frac{1}{4}\sum_{i,j,k}\widetilde\Gamma^k_{ij}e_i\cdot e_j\cdot e_k\cdot\bar\psi\;.    
\end{equation*}   
Now, using that $e_i = \sum_j b_i^j \pa_j $ and that 
$   e_i \cdot \overline{\nabla_{e_i} \psi} = \overline{\pa_i \cdot
  \nabla_{e_i} \psi}  $, we obtain that 
$$ \bar D\bar\psi=\sum_{ij} b_i^j \overline{\pa_i \nabla_{\pa_j} \psi}+
\frac{1}{4}\sum_{i,j,k}\widetilde\Gamma^k_{ij}e_i\cdot e_j\cdot
e_k\cdot\bar\psi\,$$ 
and hence,
\begin{eqnarray} \label{diracrel}
\bar D\bar\psi=\overline{D\psi}  + \sum_{ij} (b_i^j - \de_i^j) 
\overline{\pa_i \cdot  \nabla_{\pa_j} \psi}+
\frac{1}{4}\sum_{i,j,k}\widetilde\Gamma^k_{ij}e_i\cdot e_j\cdot
e_k\cdot\bar\psi.
\end{eqnarray}

\noindent Let us introduce a convenient       
notation for sections in this trivialization.       
For $v=\sum v^i e_i(P): U\to T_PM\cong\mR^n$ we define        
$\overline{v}:U\to TM$, $\overline{v}(q)=\sum v^i(q)e_i(q)$, i.e.\ $v$        
is the coordinate presentation for $\overline{v}$.        
Similarly, for $\psi:U\to \Si_PM=\Si\mR^n$,        
$\psi=\sum\psi^i\al_i(P)$ we write $\overline{\psi}(q)=\sum \psi^i(q)\al_i(q)$.       
In this notation, $\overline{x}$ is the outward radial vector field whose        
length is the radius.        
Similarly, we write $D$ for the Dirac operator on flat $\mR^n$ and        
$\overline{D}$ for the Dirac operator        
on $(M,g)$.

\begin{proposition}\label{prop.green.expand}       
Let $(M,g)$ be a compact connected spin manifold of dimension $n \in 
\{3,4,5\}$.  Let $P \in M$ and $\psi_0\in \Si_P M$, then there is a        
spinor $\Psi(\psi_0)$ which is $\bar{D}$-harmonic on $M\setminus P$, and which has   
in the trivialization defined above the following        
expansion at $P$:        
$$\Psi(\psi_0)=   \frac{\bar{ x}} {r^n} \cdot \psi_0  +         
\Th_1(x)     \hbox{ if } n=3$$        
$$\Psi(\psi_0)=          
\frac{\bar{ x}} {r^n} \cdot \psi_0 +        
\Th_2(x) \hbox{ if } n=4$$       
$$\Psi(\psi_0)=          
 \frac{\bar{ x}} {r^n} \cdot \psi_0 + \alpha(x)         
+ \Th_3(x)  \hbox{ if } n=5$$        
where $\Th_1 \in C^{0,a}(\mR^n)$ for all $a\in ]0,1[$, where, for all   
$\ep>0$,   
 $r^{\ep} \Th_2$,       
$r^{\ep} \Th_3$, $r^{1+\ep} |\na \Th_1|$, $r^{1+\ep} |\na \Th_2|$ and  $r^{1+\ep} |\na \Th_3|$ are       
continuous on  $\mR^n$ and where $\al$ is homogeneous of order $-1$ and   
even near $P$.           
\end{proposition}       
       
\begin{remark}        
As before, the spinor  $\om_{n-1}^{-1}\Psi(\psi_0)$ is the Green's        
function for the Dirac operator on $M$. The expansion of $\Psi(\psi_0)$     
could be improved but the statement of proposition \ref{prop.green.expand}       
is sufficient to adapt the proof of theorem \ref{mass}.      
\end{remark}         
        
\proof{of positive mass theorem} Using this proposition, the proof of theorem \ref{mass} can        
easily be adapted with $\psi= \Psi(\psi_0)$.        
As one can check, equation (\ref{final}) is still available in        
dimensions $3$, $4$ and $5$. This is easy to see in dimensions $3$ and   
$4$. In dimension $5$, we note that since $\al$ is even near $P$ and since      
${x\over r}$ is an odd vector field, we have      
$$Re \int_{S^{n-1}}{\pa\over \pa r}(<\frac{x}{r} \cdot \psi_0, \alpha(x) >)\, ds =         
 \int_{S^{n-1}} Re(<\frac{x}{r} \cdot \psi_0, \pa_r \alpha(x) >)\, ds=0$$        
Equation (\ref{final}) easily follows.         
This proves the positive mass theorem 2.2.        
\qed

We will now prove Proposition~\ref{prop.green.expand}.        
        
\begin{definition}       
Let $\alpha \in \Gamma(\Si(\mR^n \setminus \{0\}))$        
be a smooth spinor defined on        
$\mR^n \setminus \{0\}$. For $k \in \mR$, we say that $\alpha$ is \emph{homogeneous of order $k$}       
if $\alpha(s x)=       
s^k \alpha(x)$ for all $x\in \mR^n \setminus \{0\}$ and all $s>0$.      
This is equivalent to $\partial_r \alpha = k \alpha$.       
\end{definition}       
       
\begin{proposition} \label{phom}       
Let $\alpha$ be a spinor         
homogeneous of order $k\in (-n,-1)$. Then there is a spinor $\beta$,    
homogeneous of order $k+1$, such that $D(\beta)=\alpha$.       
\end{proposition}       
       
\proof{} Let $\alpha$ be a homogeneous spinor       
of order $k$. Recall that       
$\Gamma_D:={x \over  \om_{n-1} r^n}\cdot$ is the Green's function for the Dirac operator on  $\mR^n$.        
We define $\beta:= \Gamma_D*\al$, i.e.   
  $$\be(x)= {1 \over \om_{n-1}} \,\lim_{\rho\to 0}\int_{\mR^n\setminus (B_x(\rho)\cup B_0(\rho))}  {x-y \over |x-y|^n}\cdot \alpha(y)\, dy.\qquad x\not = 0$$   
The integral converges for $|y|\to \infty$ as $k<-1$.   
The limit for $\rho\to 0$ exists as $k>-n$.    
Similarly one sees that $\beta$ is smooth, and one calculates   
$D(\beta)=\alpha$. A simple change of variables       
shows that $\beta$ is homogeneous of order $k+1$.       
       
\begin{lemma}[Regularity Lemma] \label{phom1}       
Let $\al$ be a smooth spinor on $\mR^n\setminus\{0\}$. We       
assume that $D(\al)=O({1 \over r})$ and    
$\pa_i D(\al)=O({1 \over r^2})$ as $r\to 0$.        
Then, for all $\ep>0$, $r^\ep \al$ and $r^{1+\ep} |\na \al|$    
extend continuously to $\mR^n$.       
\end{lemma}        
       
\proof{}     
As the statement is local, we can assume for      
simplicity that $\al$ vanishes outside     
a ball $B_0(R)$ about $0$.       
Since $D(\alpha) \in L^{q}(\mR^n)$ for all $q < n$, from regularity theory we get      
that $\al \in H_1^q(\mR^n)$ for all $q < n$. The Sobolev embedding     
theorem then implies that $\al \in L^q(\mR^n)$ for all $q >1$. Moreover,      
we have      
$$| D( r^\ep \al) | \leq \ep r^{\ep-1} |\al | + O(r^{\ep-1})$$       
Using H\"older's inequality, we see that $D( r^\ep \al) \in L^{q}(\mR^n)$     
for some $q>n$. By regularity theory, we have      
$ r^\ep \al  \in  H_1^q(\mR^n)$ and by the Sobolev embedding theorem, $r^\ep \al     
\in C^0(\mR^n)$. This proves the first part of lemma \ref{phom1}. For the     
second part, we apply the same argument twice:      
a calculation on $\mR^n\setminus\{0\}$ yields:       
        
\begin{eqnarray} \label{regu0}        
|D(r^{1+\ep} \pa_i \al)| \leq  (1+\ep)  r^\ep  |\pa_i \al| +O(r^{\ep-1})        
\end{eqnarray}        
In the same way, we have:        
\begin{eqnarray} \label{regu1}        
| D( r^\ep \pa_i \al) | \leq \ep r^{\ep-1} |\pa_i \al | + O(r^{\ep-2})        
\end{eqnarray}        
For all $i=1,\ldots,n$ we have  $ D(\pa_i \al) = \pa_i D \al=  O(r^{-2})$       
and        
$D(\pa_i \al) \in L^q(\mR^n)$         
for all $q < \frac{n}{2}$. Using the regularity theory and then the Sobolev        
embedding theorem, we get that $\pa_i \al \in H_1^q(\mR^n)$ for all $q <        
\frac{n}{2}$        
and that  $\pa_i \al \in L^q(\mR^n)$ for all $q <        
n$. The H\"older inequality implies that there is a        
$q>{n\over 2}$, close  to $\frac{n}{2}$, such that $r^{\ep-1} |\pa_i \al |        
\in L^q        
(\mR^n)$. Together with (\ref{regu1}), this shows that $ D( r^\ep        
\pa_i \al)\in L^q(\mR^n)$ for some $q>\frac{n}{2}$.         
By the regularity and Sobolev theorems, we obtain        
that  $ r^\ep \pa_i \al \in L^q(\mR^n)$  for some         
$q> n$. Now using (\ref{regu0}), we obtain         
$|D(r^{1+\ep} \pa_i \al)| \in L^q(\mR^n)$ for some $q>n$.         
Applying regularity theory and the       
Sobolev theorems again, we get that $r^{1+\ep} \pa_i \al \in C^0(\mR^n)$.        
This proves that $r^{1+\ep} |\na \al| \in  C^0(\mR^n)$.       
\qed        
       
\proof{of Proposition~\ref{prop.green.expand}}       
       
\noindent 
Let us come back to formula (\ref{diracrel}). We have 
  $$\widetilde\Gamma^k_{ij}e_k=\bar\nabla_{e_i}e_j=b_i^r\bar\nabla_{\partial_r}(b_j^s\partial_s) 
    =b_i^r(\partial_r b_j^s)\partial_s+b_i^r b_j^s    
\Gamma_{rs}^l\partial_l\;,$$   
where as usually the Christoffel symbols of the first kind $\Gamma_{rs}^l$  
are defined by $$\Gamma_{rs}^l\partial_l=\bar\nabla_{\partial_r}\partial_s\;.$$   
   
\noindent Therefore we have   
$$\widetilde\Gamma^k_{ij}b_k^l\partial_l=b_i^r(\partial_r b_j^s)\partial_s+b_i^r b_j^s    
\Gamma_{rs}^l\partial_l\;,$$   
and hence   
\begin{equation}\label{gamma1}   
\widetilde\Gamma^k_{ij}=\Big(b_i^r(\partial_r b_j^l)+b_i^r b_j^s\Gamma_{rs}^l\Big)(b^{-1})^k_l\;.   
\end{equation}

\noindent Let $\eta$ be a cut-off function equal to $1$ in a        
neighborhood~$V$  
of $P=0$ in $M$, and supported in the normal neighborhood~$U$.    
Let $\psi_0$ be a constant spinor        
on $\mR^n$. We define $\psi$ on $U\setminus\{0\}$ by           
$$\overline\psi = \frac{\eta}{r^{n-1}} \frac{ \overline{x}} {r} \cdot        
\overline\psi_0        
$$        
and extend it with zero on $M\setminus U$.   
 Now,  we have the following    
development of the metric $g$ (see for example \cite{lee.parker:87}):   
   
\begin{eqnarray*}  
g_{ij} = \delta_{ij}+    
\frac{1}{3} R_{i \alpha \beta j}(p) x^{\alpha} x^{\beta}   
+O(r^3) 
 \end{eqnarray*}   
Since the matrix $(b_{ij})$ is equal to
$(g_{ij})^{-\frac{1}{2}}$, we get that 
 \begin{equation}\label{devbij}      
b_i^j=\delta_i^j-\frac{1}{6}R_{i\alpha\beta j}x^\alpha x^\beta + O(r^3).
\end{equation}
Since $Ric(p)= 0$, one computes that near $P$ (i.e. where $\eta \equiv 1$),
\begin{eqnarray} \label{or3}
\sum_{ij} R_{i\alpha\beta j}x^\alpha x^\beta  \overline{ \partial_i \cdot \nabla_{\partial j}
 \psi}= 0.
\end{eqnarray}
In the same way, using Bianchi identity and relation (\ref{gamma1}), 
we compute that 
\begin{eqnarray} \label{chr}        
\sum_{ijk}  \widetilde\Gamma_{ij}^k e_i \cdot e_j \cdot e_k  = O(r^2).       
\end{eqnarray}     
Then, $\psi$ is smooth on $M\setminus\{P\}$ and is $\bar{D}$-harmonic near $P$ (see the locally conformally flat case)        
and by  (\ref{diracrel}), near $P$ we have        
$$\overline{D} (\overline\psi) = \frac{1}{4 r^{n-1}} \sum_{ijk} \widetilde\Gamma_{ij}^k         
 e_i \cdot e_j \cdot e_k \cdot  \frac{ \overline{x}} {r} \cdot         
\overline\psi_0+  \sum_{ij} (b_i^j-\de_i^j) \overline{ \partial_i \cdot \nabla_{\partial j}
 \psi}.$$           

\noindent Writing the Taylor development of the right side member of  this
relation with the help of relations (\ref{devbij}) and(\ref{chr})  
 we can write        
$\overline{D} (\overline\psi)$ as a sum of a spinor $\gamma$ which is        
homogeneous of order $3-n$ and a spinor $\gamma'$,        
smooth on  
$U \setminus \{ 0 \}$, 
which       
satisfies $\gamma'=0(r^{4-n})$ and $|\na \gamma'|= 0(r^{3-n})$.       
  $$ \overline{D} (\overline\psi) = \gamma+ \gamma'$$        
If $n=3$, $\gamma + \gamma' \in L^q(U)$ for all $q>1$.  
Let $\eta$ be a cut-off function as above,  
then $\eta(\ga+ \ga')$ can be viewed as a spinor on  $\mR^n\setminus\{0\}$, 
and  $\Th:=\Gamma_D* (\eta(\ga+ \ga'))$ is a smooth spinor     
on $\mR^n\setminus\{0\}$ such that        
$D(\Th)= \gamma + \gamma'$ near $0$.  
By regularity theory and the Sobolev   
embedding theorem, we get that $\Th \in H_1^q(\mR^n)$ for all $q >1$,  
and hence   
$\Th \in C^{0,a}(\mR^n)$ for all $a \in (0,1)$.  
Lemma~\ref{phom1} implies $r^{1+\ep} |\na \Th|\in C^0(\mR^n)$. 
\\ 
       
\noindent If $n=4$, we have $\gamma+\gamma' = 0({1 \over r})$ and        
$|\na( \gamma + \gamma')| = 0({1 \over r^2})$. As in dimension $3$, we can   
find $\Th$, a smooth spinor     
on $\mR^n\setminus\{0\}$ such that        
$D(\Th)= \gamma + \gamma'$ near $0$.      
We fix $\epsilon \in (0,1)$.        
By the regularity lemma \ref{phom1} we get  $r^{\ep}       
\Th \in C^0(\mR^n)$ and        
$r^{1+\ep} |\na \Th| \in C^0(\mR^n)$.\\       
       
\noindent If $n=5$, by proposition \ref{phom}, we can find a        
spinor $\al$ homogeneous of order $-1$       
such that $D(\al)= \gamma$. Moreover,       
$\gamma' = 0({1 \over r})$ and $|\na\gamma'| = 0({1 \over       
  r^2})$. Proceeding as in dimension $3$ and using lemma \ref{phom1}, we can find a spinor  $f$ smooth        
on $\mR^n\setminus\{0\}$ such that        
$D(f)= \gamma'$ near $0$, such that $r^{\ep}       
f \in C^0(\mR^n)$ and such that $r^{1+\ep} |\na f| \in C^0(\mR^n)$.       
We set $\Th:=( \al + f)$.\\       
             
\noindent Now, for all dimensions, we set        
$$\varphi= \psi  - \eta \Th$$        
By (\ref{diracrel}), we have  
on $V$        
$$\overline{D} \overline\varphi =\overline{D}( \overline\psi) -    
\overline{D \Th}  -       
{1 \over 4}  \sum_{ijk}  \Ga_{ij}^k e_i \cdot e_j \cdot e_k \cdot        
\overline\Th -  \sum_{ij} (b_i^j-\de_i^j) \overline{ \partial_i \cdot \nabla_{\partial j}
 \Th}$$        
Using (\ref{or3}), (\ref{chr}) and the fact that $\overline{D}( \overline\psi) -    
\overline{D \Th}=0$, we get           
that $\overline{D} \overline\varphi =O(r)$ and hence is of class        
$C^{\infty}(M\setminus\{P\}) \cap C^{0,1}(M)$. As a consequence, there exists         
$\overline\varphi' \in \Ga(\Si M)$ of class $C^{\infty}(M\setminus\{P\}) \cap        
C^{1,a}(M)$, $a\in (0,1)$         
such that $\overline{D} \overline\varphi'= \overline{D} \overline\varphi$.    
We now        
set $\Psi(\psi_0)=\overline{\varphi}-  \overline\varphi'$.        
Proposition~\ref{prop.green.expand} follows.\qed

\vspace{1cm}           
Authors' address:           
\vspace{5mm}           
\parskip0ex           
{\obeylines           
Bernd Ammann and Emmanuel Humbert,           
Institut \'Elie Cartan BP 239          
Universit\'e de Nancy 1           
54506 Vandoeuvre-l\`es -Nancy Cedex           
France           
}           
\vspace{0.5cm}           
           
E-Mail:           
{\tt ammann at iecn.u-nancy.fr, humbert at iecn.u-nancy.fr}           
   
\end{document}